\documentclass[12pt]{amsart}
\usepackage{amscd}
\newtheorem{thm}{Theorem}
\newtheorem{prop}{Proposition}
\newtheorem{lem}{Lemma}
\newtheorem{cor}{Corollary}

\theoremstyle{remark}
\newtheorem{rem}{Remark}

\theoremstyle{definition}

\newtheorem{ex}{Example}

\newcommand{\GG}{\mathcal{ G}}

\newcommand{\MA}{\mathcal{ A}}
\newcommand{\BB}{\mathcal{ B}}

\newcommand{\Q}{\mathbb{ Q}}

\newcommand{\N}{\mathbb{ N}}

\newcommand{\rk}{\operatorname{rk}}

\newcommand{\Z}{\mathbb{ Z}}

\title[Rational topology of gauge groups]{The rational topology of gauge groups and of spaces of connections}
\author{Svjetlana Terzi\'c}
\address{Faculty of Science, University of Montenegro,
Cetinjski put bb, 81000 Podgorica, Montenegro, Serbia and Montenegro
}
\email{sterzic@cg.ac.yu}
\thanks{The results in this paper are partially obtained while the author was postdoc supported by the
{\it DFG Graduiertenkolleg ``Mathematik im Bereich
ihrer Wechselwirkung mit der Physik''} at the Mathematical Institute of the
Ludwig-Maximilians University in Munich}
\date{\today; MSC 55P62, 57R19, 58B05, 81T13}

\begin{document}

\begin{abstract}
Let $P$ be a principal bundle with semisimple compact simply connected
structure group $G$ over a compact simply connected four-manifold $M$.
In this note we give  explicit formulas for the rational
homotopy groups and cohomology
algebra of the gauge group and of the space of (irreducible) connections modulo gauge
transformations for any such bundle.
\end{abstract}

\maketitle
\begin{center}
{\bf To the memory of my teacher and friend \\ Professor Yuri Petrovich
Solovyov}
\end{center}

\section{Introduction}\label{s:intro}
Our aim in this note is to compute the rational homotopy groups and
the rational cohomology of the gauge group and of the space of connections modulo gauge
transformations for principal bundles over four-manifolds.
We assume that $G$ is a semisimple
compact simply connected Lie group and $M$ is  compact and simply connected.

Note that in some particular cases some  of these computations have already
been done. Namely, Donaldson~\cite{don}, whose original calculation proceeds the book~\cite{DK},
computed the cohomology structure of the quotients of spaces of connections for
SU(2)-principal bundles over compact simply connected four-manifolds,
but the proof essentially uses the fact that the structure group is $SU(2)$.

We propose here a general approach which appeals to Sullivan's minimal model
theory. We proceed as follows.
First we compute the rational homotopy groups of the
gauge group using the result of Singer~\cite{singer} characterising
the weak
homotopy type of base point preserving gauge groups,
and Whitehead-Milnor theorem~\cite{M},~\cite{W} giving the homotopy type of
a simply connected four-manifold $M$. Then the corresponding fibrations
between the gauge groups and spaces of (irreducible) connections
yield the  rational homotopy groups of the quotients.

Having computed the rational homotopy groups of the quotients
of spaces of connections,
the nilpotency of the space of
connections modulo based point gauge transformations group will make it
possible to apply Sullivan's minimal model theory for the
cohomology computation.
Since it will turn out that these spaces have free cohomology algebras, it
will immediately imply that they are formal in the sense of Sullivan.

In the case of the gauge groups, the application of the Sullivan minimal
model theory is possible because they are  $H$-spaces.

As an examples we recover the calculation of the
rational cohomology  in the case
when $G = SU(2)$ and also provide the calculation in the case when $G=E_8$.

\subsection{General remarks.}
In this section we briefly recall some general facts on the topology
of gauge groups and spaces of connections. For definitions and detailed
information on this topic we refer to~\cite{DK}.

Let $\pi : P\to M$ be a $G$ - principal bundle, where $G$ is
a semisimple compact simply connected Lie group and $M$ a compact simply
connected four-manifold.    Let $\GG$ be the group of gauge
transformations
of this principal bundle which induce the  identity on the base.

Let $\MA$ denote the space of all connections and $\MA ^{*}$ the subspace
 of all
irreducible connections on the bundle $P$. We assume that $\MA $
and $\GG$ are
equipped with certain Sobolev topologies.
Usually one fixes the Sobolev class  $L_{p-1}^{2}$ for $\MA$ and
$L_{p}^{2}$ for $\GG$, but since we
will be interested only in homotopy invariant properties, the  particular
choice of $p$ is not important, as long as $p$ is large enough.

The action of $\GG$ on $\MA$ and $\MA ^{*}$ is not free in general.
 In order to have a free action on $\MA$, it turns out
that one should consider those gauge transformations which fix one fiber.
Such automorphisms of $P$ we denote by $\GG _{0}$.
Besides that, one also has the  free action of the group
$\tilde {\GG} = \GG /Z(G)$ on $\MA ^{*}$, where by $Z(G)$ we denote the
center of the group $G$.
Note that $Z(G)$ is finite, since we assume $G$ to be semisimple.

We denote by $\tilde \BB = \MA /\GG _{0}$, then
$\tilde{\BB}^{*} = \MA ^{*}/\GG _{0}$ and $\BB ^{*} = \MA ^{*}/\tilde {\GG}$
the corresponding quotients.
Since in all these cases we have free actions, one gets the corresponding
fibrations
\begin{equation}\label{TB}
\GG _{0}\to \MA \to \tilde \BB \ ,
\end{equation}
\begin{equation}\label{TB*}
\GG _{0}\to \MA ^{*}\to \tilde{\BB}^{*} \ ,
\end{equation}
\begin{equation}\label{B*}
\tilde{\GG}\to \MA ^{*}\to \BB ^{*} \ .
\end{equation}

We are interested in the rational topology of these objects, i. e. firstly
in a computation of their rational homotopy groups, and secondly in a
computation of their rational cohomology.

For these  purposes the following observations are useful, see~\cite{DK}.
\subsection{Remarks}\label{obs}
\begin{itemize}
\item $\MA$ is contractible, since it is an affine space;
\item $\pi _{j}(\MA ^{*}) = 0$ as $\MA ^{*}$ is weakly homotopy equivalent to
$\MA$;
\item ${\tilde \BB}^{*}$ is weakly homotopy equivalent to ${\tilde \BB}$;
\item $\tilde{\GG} = \GG /Z(G)$ and hence we have a fibration
$Z(G)\to \GG \to \tilde{\GG}$ which gives an exact homotopy sequence
$$
...\to \pi _{j}(Z(G))\to \pi _{j}(\GG )\to \pi _{j}(\tilde{\GG})\to
\pi _{j-1}(Z(G))\to ...
$$
Since $Z(G)$ is finite we  have that
$
\pi _{j}(Z(G)) = 0 \; \mbox{for} \; j\geq 1 \
$
and, therefore
$$
\pi _{j}(\GG ) = \pi _{j}(\tilde{\GG}) \;  \mbox{for} \;  j\geq 2 \ .
$$

For $j=1$ we have the exact sequence
\begin{equation}\label{first}
0\to \pi _{1}(\GG)\to \pi _{1}(\tilde{\GG})\to Z(G)\to \pi _{0}(\GG )\to
\pi _{0}(\tilde{\GG}) \ .
\end{equation}

\item  $G = \GG /\GG _{0}$ and we have  also the fibration
 $\GG _{0}\to \GG \to G$,  which implies the exact
sequence
\begin{equation}\label{second}
\pi _{j+1}(G)\to \pi _{j}(\GG _{0})\to \pi _{j}(\GG )\to \pi _{j}(G) \ .
\end{equation}
\end{itemize}

\subsection{ Algebraic topology tools.}~\label{att}
We recall here some well known facts from algebraic topology which will
be useful for our purposes,
see~\cite{switzer}.

Let $M = X\cup _{f}CY$ be a cofibration, i.e. intuitively we glue the
base of the cone $CY$ to $X$ by means of $f$. Denote by $SX$ the suspension
of $X$ which we may regard as being the quotient
$(I\times X)/(\{ 0\}\times X\cup I\times \{ x_0 \}\cup \{ 1 \}\times X).$
Then the sequence
$$
(Y, *)\stackrel{f}\to
 (X, *)\stackrel{j}\to (M, *)\stackrel {k^{'}}\to (SY, *)
\stackrel{Sf}\to (SX, *)\stackrel{Sj}\to (SM, *)\to \ldots
$$
$$
\ldots \stackrel{S^{n}f}\to (S^nY, *)\stackrel{S^{n}f}\to (S^{n}X, *)
\stackrel{S^{n}j}\to (S^{n}M, *)\to \ldots
$$
is coexact. This means that the  mapping functor $Map(-, G)$ will turn this
sequence  into an exact sequence
\begin{equation}\label{cofib}
[Y;G]\stackrel{f^{*}}\leftarrow [X;G]\stackrel{j^{*}}\leftarrow
 [M;G]\stackrel{k^{'*}}
\leftarrow [SY;G]\stackrel{Sf^{*}}\leftarrow [SX;G]
\stackrel{Sj^{*}}\leftarrow [SM;G]\leftarrow \ldots
\end{equation}

The following obvious observation will be useful for our further
computation.

\begin{lem}\label{lem}
\begin{equation}
\pi _{j}(Map_{*}(M, G)) = [S^jM; G], j\in {\mathbb N} \ .
\end{equation}
\end{lem}

\section{Rational homotopy groups of gauge groups and of spaces of
 connections}

Let  $M$ be a compact simply connected four-dimensional manifold.
Then, by the result of Whiethead-Milnor~\cite{M},~\cite{W}, we know that $M$ is homotopically
a cofibration, i. e. $M=\vee _{b_2(M)}S^2\cup _{h}D^4$.
 It means that, up to homotopy, $M$ is obtained by attaching a four-cell
to the wedge of two-spheres by an attaching map
 $h : S^{3}\to \vee _{b_2(M)}S^2$. Since $D^{4}=CS^{3}$ and $SS^{k}=S^{k+1}$,
we can
apply what was said in~\ref{att} to this
cofibration and get the following exact sequence
$$
[S^3;G]\leftarrow [\vee _{b_2(M)}S^2;G]\leftarrow [M;G]\leftarrow [S^4;G]
\leftarrow [\vee _{b_2(M)}S^3;G]\leftarrow \ldots
$$
$$
\leftarrow [S^{n+3};G]
\leftarrow [\vee _{b_{2}(M)}S^{n+2};G]\leftarrow [S^nM;G]\leftarrow \ldots \ ,
$$
which is actually
\begin{equation}\label{bs}
\pi _{3}(G)\leftarrow \oplus _{b_2(M)}\pi _{2}(G)\leftarrow [M;G]
\leftarrow \pi _{4}(G)\leftarrow \oplus _{b_2(M)}\pi _{3}(G)\leftarrow
\ldots
\end{equation}
$$
\leftarrow  \pi _{n+3}(G)\leftarrow \oplus _{b_2(M)}\pi _{n+2}(G)\leftarrow
 [S^{n}M;G]\leftarrow \ldots
$$

We are  able to compute first the rational homotopy groups of $\GG _{0}$
because of the following theorem of Singer~\cite{singer} which,
among  other things,
proves that the weak homotopy type of $\GG _{0}$ is independent of $P$.

\begin{thm}\label{singer}
If $G$ is a compact simply connected semisimple Lie group, then for
$\dim M\leq 4$ we have that
\begin{equation}
\GG _{0}\sim Map _{*}(M, G) \ ,
\end{equation}
i.e.
$\GG _{0}$ is weakly homotopic to the space $Map_{*}(M, G)$.
\end{thm}

In particular,  Theorem~\ref{singer} gives that in the case of
four-dimensional manifolds we have that
\[
\pi _{j}(\GG _{0}) = \pi _{j}(Map _{*}(M, G)) \ ,
\]
for any $j\in \mathbb{N}$.

Let us first recall some well known facts on the rational
cohomology of connected compact
Lie groups.

\begin{rem}
The classical Hopf theorem~\cite{borel} gives that the cohomology
algebra of a compact connected
Lie group is an  exterior algebra over odd degree generators, i. e.
$H^{*}(G)=\wedge (z_1,\ldots ,z_n)$, where $n=\rk G$ and $\deg z_i = 2k_i-1$.
Here $\rk G$ is the rank of the group $G$, i.~e. the dimension of the maximal torus in $G$ and
the numbers $k_i$, $1\leq i\leq n$ are so called exponents of the group $G$.
Because of formality,  rational homotopy theory gives then that
 $\rk \pi _{2k}(G)=0$, that $\rk \pi _{2k-1}(G)=0$ if $k$ is not an exponent
for $G$, while for a $k$ being an exponent for $G$ we have that
$\rk \pi _{2k-1}(G) = b_{2k-1}(G)= \nu (k)$. Here by $\nu (k)$ we denote
the multiplicity of the exponent $k$. Obviously, it is satisfied
that $\rk G=\sum_{j\in \N}\rk \pi _{j}(G)$.
\end{rem}

\begin{rem}
Since $\pi _{2j}(G)$ are finite, Theorem~\ref{singer} and
the sequence~\eqref{bs} implies that $\pi _{0}(\GG _0)$ is finite.
\end{rem}

\begin{rem}
Since we are interested in the rational homotopy groups of these spaces we
need to  consider the sequence~\eqref{bs} tensored by $\Q$. Since
$\pi _{1}(\GG _{0})$ is abelian, we can also tensor
$\pi _{1}(\GG_{0})$ by $\Q$.
Namely, $\GG _{0}$ is topological group and, thus,  all its
connected components are homeomorphic, so we can fix the  component
$\GG _{0}^{e}$
corresponding to the identity automorphism. It is obviously also
topological group and, thus, homotopy simple.\footnote{Meaning that its
fundamental group is abelian and acts trivially on its higher
homotopy groups.}
%Then
%$\pi _{1}(\GG _{0})=\oplus _{\vert \pi _{0}(\GG )\vert}\pi_{1}(\GG _{0}^{e})$.
%Because of Theorem~\ref{singer} we have that
%$\GG _{0}^{e}$ is weakly homotopy
%equivalent to
%the corresponding connected component in $Map _{*}(M, G)$.
%Now, since $M$ and $G$
%are connected,  the result of Hilton~\cite{hilton} implies that any connected
%component of $Map _{*}(M, G)$ is nilpotent.
%Because of the weak homotopy
%equivalence, this gives that $\GG _{0}^{e}$ is nilpotent.
\end{rem}

\begin{prop}\label{rgo}
The ranks of the homotopy groups of the group $\GG _{0}$ are given by

\begin{equation}
\rk \pi _{j}(\GG _{0}) = b_2(M)\rk \pi _{j+2}(G) +
\rk \pi _{j+4}(G), \;  j\in \mathbb{N} \ .
\end{equation}

\end{prop}
\begin{proof}

Applying Theorem~\ref{singer} to the sequence~\eqref{bs} we have
that
\begin{equation}\label{seq}
\ldots \leftarrow
\pi _{j+3}(G)\otimes \Q \leftarrow \oplus _{b_2(M)}\pi _{j+2}(G)\otimes \Q
\leftarrow
\pi _{j}(\GG _{0})\otimes \Q \leftarrow \pi _{j+4}(G)\otimes \Q \ldots
\end{equation}

Since for $j$ even we know that  $\pi _{j}(G)\otimes \Q = 0$, the
sequence~\eqref{seq} immediately implies
that $\pi _{j}(\GG _{0})\otimes \Q =0$ for $j$ even.

Therefore, we will assume that $j$ is odd.
Then from the
sequence~\eqref{seq} we can extract the following short exact sequence

$$
0\leftarrow b_{2}(M)\pi _{j+2}(G)\otimes \Q \leftarrow \pi _{j}(\GG _0)\otimes \Q
\leftarrow \pi _{j+4}(G)\otimes \Q \leftarrow 0 \ .
$$
This gives that
\begin{equation}
\rk \pi _{j}(\GG _{0}) = b_2(M)\rk \pi _{j+2}(G) +
\rk \pi _{j+4} (G) \ .
\end{equation}
\end{proof}

\begin{prop}\label{rg}
The ranks of the homotopy groups of the group $\GG$ are given by
\begin{equation}
\rk \pi _{j}(\GG ) = b_2(M)\rk \pi _{j+2}(G) + \rk \pi _{j+4}(G) +
\rk \pi _{j}(G), j\in \mathbb{N} \ .
\end{equation}
\end{prop}

\begin{proof}
From the same reasons as in the case of the group $\GG _{0}$
we are able to tensor the sequence~\eqref{second} by $\Q$.

Thus, tensoring the  sequence~\eqref{second} by $\Q$ we get the exact sequence
\begin{equation}\label{zag}
\pi _{j+1}(G)\otimes \Q \to \pi _{j}(\GG _{0})\otimes \Q \to
\pi _{j}(\GG )\otimes \Q \to \pi _{j}(G)\otimes \Q \to
\pi _{j-1}(\GG _0)\otimes \Q \ .
\end{equation}

First, for $j$ even this sequence implies that
$\pi _{j}(\GG)\otimes \Q =0$.

Therefore we assume $j$ to be odd. Then the sequence~\eqref{zag} and
Proposition~\ref{rgo} give the following short exact sequence
$$
0\to \pi _{j}(\GG _0)\otimes \Q \to \pi _{j}(\GG)\otimes \Q \to
\pi _{j}(G)\otimes \Q \to 0 \ .
$$
This  implies that
\begin{equation}
\rk \pi _{j}(\GG ) = \rk \pi _{j}(G) +
\rk \pi _{j}(\GG _0)
\end{equation}
$$
= b_{2}(M)\rk \pi _{j+2}(G) + \rk \pi _{j+4}(G)
+ \rk \pi _{j}(G) \ .
$$

\end{proof}

\begin{rem}\label{rt}
Since we assume $G$ to be semisimple, the
sequences~\eqref{first} and~\eqref{second} give that for all $j\in \mathbb{N}$
$$
\pi _{j}(\tilde{\GG})\otimes \Q = \pi _{j}(\GG)\otimes \Q \ ,
$$
and by the above Proposition we have computed
 the rational homotopy groups for
$\tilde {\GG}$.
\end{rem}

\begin{cor}\label{rb}
The rational homotopy groups for ${\tilde {\BB}}$, ${\tilde {\BB}}^{*}$   and
$\BB ^{*}$ are given by
\begin{itemize}
\item $\rk \pi _{j}({\tilde {\BB}})
      = \rk \pi _{j}({\tilde {\BB}}^{*}) = b_{2}(M)\rk \pi _{j+1}(G)
   +  \rk \pi _{j+3}(G)$, $j\geq 1$ \ ,
\item $\rk \pi _{j}(\BB ^{*}) = b_{2}(M)\rk \pi _{j+1}(G)
   + \rk \pi _{j+3}(G) + \rk \pi _{j-1}(G)$, $j\geq 1$ \ .
\end{itemize}
\end{cor}

\begin{proof}
Since the homotopy groups of the
total spaces in the fibrations~\eqref{TB},~\eqref{TB*}
and~\eqref{B*} are trivial we get that
$$
\pi _{j}({\tilde {\BB}})\otimes \Q = \pi _{j}({\tilde {\BB}}^{*})\otimes \Q =
\pi _{j-1}(\GG _{0})\otimes \Q \ ,
$$
$$
\pi _{j}(\BB ^{*})\otimes \Q = \pi _{j-1}({\tilde {\GG}})\otimes \Q \ .
$$
Now, the Propositions~\ref{rgo},~\ref{rg} and Remark~\ref{rt}
give the statement.
\end{proof}

\section{Rational cohomology of gauge groups and of spaces of connections}

\subsection{Rational cohomology of gauge groups.}

As we already said in the proof of Proposition~\ref{rg}, the identity
component $\GG ^{e}$ of the gauge group $\GG$  is homotopy simple space.
By the Proposition~\ref{rg} it is of finite type and has only odd
degree nontrivial rational homotopy groups.

Now we can apply Sullivan's minimal model theory, since, more generally, for
nilpotent spaces\footnote{They are given by the condition that the fundamental group is nilpotent and acts nilpotently
on higher homotopy groups.} of finite type it works well.
Namely, for a nilpotent space $X$
of finite type minimal model theory gives that the degrees and the numbers
of generators in its minimal model are given by its nontrivial
rational homotopy groups, see~\cite{L}.
More precisely, if $\mu (X)$ denotes the minimal model for $X$ then the
 number of its generators $\mu _{j}(X)$ of degree $j$ is equal to
$\rk \pi _{j}(X)$.

In our case this gives that the minimal model for
$\GG ^{e}$ has only odd degree generators.
On the other hand, for $H$-spaces of finite type the Hopf theorem~\cite{borel}
implies that their cohomology algebra is a free commutative
algebra. Therefore, $H^{*}(\GG ^{e})$
is an exterior  algebra of odd degree generators.
In particular this gives that $\GG ^{e}$ is formal in the sense of Sullivan and its minimal
model coincides with its cohomology algebra.
Moreover, $H^{*}(\GG ^{e} )$
has generators of degree $2j-1$ if and only if $\rk \pi _{2j-1}(\GG )\neq 0$.

Obviously, the cohomology algebra $H^{*}(\GG )$ is equal to the sum of
$\vert \pi _{0}(\GG )\vert$ copies of $\GG ^{e}$, i.~e.
\begin{equation}\label{cg}
H^{*}(\GG ) = \oplus _{|\pi _{0}(\GG )|}H^{*}(\GG ^{e}) \ .
\end{equation}

\begin{rem}\label{pov}
If for the structure group $G$ we have that $\pi _{4}(G)=0$, then
the sequence~\eqref{bs} implies that $\GG _{0}$ is connected. That
is the case, say, when $G = SU(n)$, $n\geq 3$ $G = Spin(n)$, $n\geq 6$
or $G$ is simply connected Lie group of exceptional type.
Obviously, the same is true for the groups $\GG$ and $\GG _{0}$.
%We know that for a topological group the rank of  any homotopy group
%has to be a multiple of the number of its connected components.
%If we apply this
%to the group $\GG$ we conclude the following.
%If the largest exponent $j_0$ of the structure group $G$ is of multiplicity
%$1$ (which is the case for simple compact Lie groups) Proposition~\ref{rg}
%gives that $\rk \pi _{2j_0 -1}(\GG ) =1$. This implies that in this case
%$\GG $ has to be connected, and, thus nilpotent space.
\end{rem}

The equation~\eqref{cg} and Proposition~\ref{rg} immediately
give the following.

\begin{thm}\label{tcg}
If $G$ is a compact simply connected semisimple Lie group,
then $H^{*}(\GG ^{e} )$
is an exterior algebra
in $(b_2(M)+2)\rk G -1$ odd degree generators. The number of generators of
degree $j$ is equal to
$b_2(M)\rk \pi_{j+2}(G)+\rk \pi_{j+4}(G)+\rk \pi_{j}(G)$.
\end{thm}

%\begin{rem}
%Being  quotients of $\GG$, under the assumptions of the Remark~\cite{pov},
%the groups $\GG _{0}$ and ${\tilde {\GG}}$ are also connected.
%\end{rem}

\subsection{Rational cohomology of ${\tilde {\BB}}$ and of
${\tilde {\BB}^{*}}$}

\begin{thm}\label{cbt}
The rational cohomology algebra of ${\tilde {\BB}}$ is a polynomial algebra in
 $(b_{2}(M)+1)\rk G - 1$ generators of even degree.
The number of generators of
degree $j$ is equal to $b_2(M)\rk \pi_{j+1}(G) + \rk \pi_{j+3}(G)$.
\end{thm}

We give a general proof which works for any semisimple compact
simply connected Lie group $G$. As we point out in the Remark~\ref{cs} below,
under some additional assumptions on the group $G$, the proof is much simpler.

\begin{proof}[Proof of  Theorem~\ref{cbt}]
By the result of~\cite{atbott} and~\cite{don}, we know
that $\tilde {\BB}$ has the weak homotopy type of
$Map_{*}(M, B_{G})_{P}$, where the latter denotes the homotopy class
of $Map_{*}(M, B_{G})$ corresponding to the maps inducing the bundle $P$.
Since $M$ and $B_{G}$ are connected,  the result of
Hilton~\cite{hilton}
implies that  any connected component
of the space $Map_{*}(M, B_{G})$ is nilpotent. Hence, $Map_{*}(X, B_{G})_{P}$
is also nilpotent.
As $\tilde {\BB}$ is weakly homotopy equivalent to
$Map _{*}(X, B_{G})_{P}$, it implies that
  $\tilde {\BB}$ is also nilpotent.
Moreover, by Corollary~\ref{rb} we know
that ${\tilde \BB}$ is of finite type, so we can apply Sullivan's
minimal model theory.
It gives that for all $j$
$$
\rk \pi _{j}({\tilde {\BB}})=\dim (\mu ({\tilde {\BB}}))_{j} \ ,
$$
and, thus
$$
\mu ({\tilde \BB})= \Q [x_1,\ldots ,x_p] \ ,
$$
where $x_i$ correspond to the non-trivial rational homotopy groups of
${\tilde {\BB}}$. Since, by  Corollary~\ref{rb},
 $\pi _{j}({\tilde {\BB}})\otimes \Q = 0$ for $j$ odd,
we have that all $x_i$ are of even degree.
This implies  that the differential in this minimal
model has to be zero.
Since by the definition of the minimal model
we have that
$H^{*}(\mu ({\tilde {\BB}}))\cong H^{*}({\tilde {\BB}})$, it follows that
\begin{equation}
H^{*}(\tilde {\BB},\Q )=\Q [x_1,\ldots ,x_p] \ .
\end{equation}
By Corollary~\ref{rb} the number of generators of degree $j$ ($j\geq 1$) is\\
$b_{2}(M)\rk \pi_{j+1}(G)+\rk \pi_{j+3}(G)$.
Then, if we sum these numbers, we get
that the number of generators in $H^{*}({\tilde \BB}, \Q )$ is
$(b_{2}(M)+1)\rk G -1$.
\end{proof}

Since the above Theorem gives that $H^{*}({\tilde {\BB}})$ is a free algebra,
it immediately implies the following.
\begin{cor}
The spaces $\tilde \BB$ and ${\tilde \BB}^{*}$ are formal in the sense of Sullivan.
\end{cor}

\begin{rem}
We pointed out in~\ref{obs} that ${\tilde {\BB}^{*}}$ is
weakly homotopy equivalent to
${\tilde {\BB}}$, hence  they have the same rational cohomology.
\end{rem}

\begin{rem}\label{cs}
If for the structure group $G$ we have that $\pi _{4}(G)=0$, then from
Remark~\ref{pov} and the fibrations~\eqref{TB} and~\eqref{TB*} it follows that
the spaces ${\tilde {\BB}}$ and ${\tilde {\BB}}^{*}$ are simply connected.
Therefore it follows immediately that the generators in the minimal models
for ${\tilde {\BB}}$ and ${\tilde {\BB}}^{*}$ are given by their nontrivial
rational homotopy groups.
\end{rem}

\subsection{Rational cohomology of $\BB ^{*}$.}

\begin{thm}\label{cb*}
The rational cohomology algebra of $\BB ^{*}$ is a polynomial algebra
in even degree $(b_2(M)+2)\rk G -1$ generators. The number of generators of
degree $j$ is equal $b_2(M)\rk \pi_{j+1}(G)+\rk \pi_{j+3}(G)+\rk \pi_{j-1}(G)$.
\end{thm}

\begin{proof}
We can get $\BB ^{*}$ as a base of the principal fibration
${\tilde \BB}^{*}\to  {\BB}^{*}$ with
a fiber $G/Z(G)$. Note that, since $Z(G)$ is finite, the groups $G$ and
$G/Z(G)$, and hence their classifying spaces,
have the same rational cohomology.
Therefore, without loss of generality, we are not going to
differentiate between them.

Since we are interested in the cohomology of ${\BB}^{*}$,
it turns out to be much easier to apply the
Borel construction to this fibration
in order to work with the total space of the fibration instead of the base.
This goes as follows.
Let $E_{G}\to B_{G}$ be the universal fibration for the
group $G$ and
denote
by ${\tilde \BB}_{G}^{*}$ the quotient
${\tilde \BB}_{G}^{*}=({\tilde \BB}^{*}\times E_{G})/G$ given
by the diagonal action of the group $G$.  Let us consider
the Serre fibration
$$
{\tilde \BB}_{G}^{*}\to B_{G} \ ,
$$
 with a fiber ${\tilde \BB}^{*}$. Obviously ${\tilde \BB}_{G}^{*}$ is
 weakly homotopy equivalent to $\BB ^{*}$.

Let $\mu ({\tilde \BB}^{*})$ and $\mu (B_{G})$ be the
minimal models for ${\tilde \BB}^{*}$ and $B_{G}$ respectively. The fact
 that $B_{G}$ is  simply connected implies that
$(\mu ({\tilde \BB}^{*})\otimes \mu (B_{G}), d)$ is a Sullivan model for
${\tilde \BB}_{G}^{*}$, see~\cite{FHT}. Moreover,  since all the generators
in $\mu ({\tilde \BB}^{*})$ and $\mu (B_{G})$ are of even degree, it follows
that $d=0$ and $(\mu ({\tilde \BB}^{*})\otimes \mu (B_{G}), d = 0)$
is the minimal model for
${\tilde \BB}^{*}_{G}$ and, hence, for  $\BB ^{*}$. It implies that
$$
H^{*}(\BB ^{*}) = H^{*}({\tilde {\BB}})\otimes H^{*}(B_{G}) \ .
$$
It is a well known fact~\cite{borel}
 that $H^{*}(B_G)$ is a polynomial algebra in
$\rk G$ generators
and the number of its generators of degree $j$ is equal $\rk \pi _{j}(G)$.
Combining this with  Theorem~\ref{cbt}, the statement follows.
\end{proof}
Since the previous Theorem gives that $H^{*}(\BB ^{*})$ is a free algebra,
it implies
\begin{cor}
The space $\BB ^{*}$ is formal in the sense of Sullivan.
\end{cor}

\begin{rem}
Again, if for the structure group $G$ we have that $\pi _{4}(G)=0$, then $\BB ^{*}$ is simply connected
and using Proposition~\ref{rb}  we immediately get its minimal model.
\end{rem}

As a direct application of the above results
we will get the rational cohomology of the gauge group  $\GG$ and
the Proposition proved in~\cite{DK} on rational cohomology
of ${\tilde \BB}$ and $\BB ^{*}$ for $SU(2)$ principal bundles over
four-manifolds. We also provide the same calculation when
the structure group of the bundle is $E_8$.

\begin{ex}
Let $G=SU(2)$. It has one exponent $k=2$ of multiplicity $1$ and, thus,
$\pi_{3}(G)\otimes \Q =\Q$,
while all the other rational homotopy groups are trivial.
Proposition~\ref{rgo} gives that
$$
\pi _{j}(\GG _0)\otimes \Q =0,\;\mbox{for}\; j\geq 2,\;
 \pi _{1}(\GG _0)\otimes \Q =\Q ^{b_2(M)} \ .
$$
Further by Proposition~\ref{rg}, we get now that
$$
\pi _{1}(\GG )\otimes \Q =\Q ^{b_2(M)} ,\;\mbox \;\pi _{3}(\GG )\otimes \Q =\Q \ ,
$$
while all the others rational homotopy groups for $\GG$ are trivial. The same
is true for $\tilde{\GG}$.

Using Corollary~\ref{rb}
we get that the non-trivial rational homotopy groups for $\tilde {\BB}$
and $\BB ^{*}$
are given by
$$
\pi _{2}(\tilde {\BB})\otimes \Q = \Q ^{b_2(M)}  \ ,
$$
$$
\pi _{2}(\BB ^{*})\otimes \Q =\Q ^{b_2(M)} ,\;\;
\pi _{4}(\BB^{*})\otimes \Q=\Q \ .
$$
Then Theorem~\ref{cbt} says that the cohomology algebra for
${\tilde \BB}$ is given by
$$
H^{*}(\tilde {\BB}) = \Q [x_1,\ldots ,x_{b_2(M)}], \;\; \deg x_i =2 \ .
$$
Theorem~\ref{cb*} gives that the cohomology algebra for $\BB ^{*}$
is
$$
H^{*}(\BB ^{*}) = \Q [x_1,\ldots ,x_{b_2(M)}, y], \;\; \deg x_i =2, \;
\deg y = 4 \ .
$$
Also, by Theorem~\ref{tcg} we have that
$$
H^{*}(\GG ^{e}) = \wedge (z_1,\ldots ,z_{b_2(M)}, w), \;\; \deg z_i = 1, \;
\deg w =3 \ .
$$

\end{ex}

\begin{ex}
It is also interesting  to consider the case when $G$ is simply connected exceptional Lie group
of type $E_8$.
The  non-trivial rational homotopy groups
for $E_8$ are given by $\pi _{2k-1}(E_8)\otimes \Q = \Q$
where $k=2,8,12,14,18,20,24,30$, since these are the exponents for $E_8$.
Proposition~\ref{rgo} gives that the non-trivial rational
homotopy groups for $\GG _{0}$ are given by
$$
\pi _{j}(\GG _{0})\otimes \Q = \Q,\; \mbox{for} \; j=11, 19, 23, 31, 35, 43, 55;
$$
$$
\pi _{j}(\GG _{0})\otimes \Q = \Q ^{b_{2}(M)},\; \mbox{for} \; j=1, 13, 21, 25,
33, 37, 45, 57 \ .
$$
Then Proposition~\ref{rg} gives that the non-trivial rational homotopy
groups for $\GG$ (and for $\tilde {\GG}$) are
$$
\pi _{j}(\GG)\otimes \Q = \Q,\; \mbox{for} \; j=3, 11, 15, 19, 27, 31,
39, 43, 47, 55, 59;
$$
$$
\pi _{j}(\GG)\otimes \Q = \Q \oplus \Q,\; \mbox{for}\; j=23, 35;
$$
$$
\pi _{j}(\GG)\otimes \Q = \Q ^{b_{2}(M)},\; \mbox{for}\; j= 1, 13, 21, 25, 33,
37, 45, 57 \ .
$$

By Corollary~\ref{rb} we see that the non-trivial rational
homotopy groups for $\tilde {\BB}$ and $\BB ^{*}$ are given by
$$
\pi _{j}(\tilde {\BB})\otimes \Q = \Q,\; \mbox{for}\; j=12, 20, 24, 32,
36, 44, 56;
$$
$$
\pi _{j}(\tilde {\BB})\otimes \Q = \Q ^{b_{2}(M)},\; \mbox{for}\; j = 2, 14,
22, 26, 34, 38, 46, 58;
$$
$$
\pi _{j}(\BB ^{*})\otimes \Q = \Q,\;
\mbox{for}\; j = 4, 12, 16, 20, 28, 32, 40, 44, 48, 56, 60;
$$
$$
\pi _{j}(\BB ^{*})\otimes \Q = \Q \oplus \Q ,\;
\mbox{for}\; j = 24, 36;
$$
$$
\pi _{j}(\BB ^{*})\otimes \Q = \Q ^{b_{2}(M)},\; \mbox{for}\; j = 2, 14, 22,
26, 34, 38, 46, 58 \ .
$$
Then by Theorem~\ref{cbt}, the cohomology algebra
$H^{*}(\tilde {\BB})$ is polynomial algebra in $8b_{2}(M)+7$ generators whose degrees
are given by the non-trivial rational homotopy groups for $\tilde
{\BB}$. Analogously, by Theorem~\ref{cb*} the cohomology algebra
$H^{*}(\BB ^{*})$ is polynomial algebra in $8b_{2}(M)+15$
generators whose degrees are given by the non-trivial rational
homotopy groups for $\BB ^{*}$. Finally, Theorem~\ref{tcg} gives the
cohomology algebra $H^{*}(\GG ^{e})$ is an exterior algebra
in $8b_{2}(M)+15$ generators whose degrees are also given by the
non-trivial rational homotopy groups for $\GG$.
\end{ex}

\begin{rem}
Note that the results of this paper do not hold in general if one omits  the assumption that the
structure group is simply connected. Namely, let us, for example, take $P$ to be
a $U(3)$ - principal bundle over $S^4$. According to~\cite{atbott}, $\tilde {\BB}$ has weak homotopy type of
$K(\Z ; 2)\times K(\Z ; 2)\times K(Z ; 4)\times K(\Z ; 6)$.
Therefore, the rational cohomology algebra for $\tilde {\BB}$ is given by
$$
H^{*}(\tilde {\BB}) = \Q [x_1, x_2, y, z ]\ ,
$$
where $\deg x_1 =\deg x_2 =2, \deg y = 4$  and $\deg z= 6$.

On the other hand, Theorem~\ref{cbt} would give that the number of generators in
$H^{*}(\tilde{\BB})$ is $2$, and, thus, it is not valid in this case.
This also shows that in the Theorem~\ref{singer} of Singer,  the condition that the group $G$ is
simply connected can not be omitted.

\end{rem}

\bibliographystyle{amsplain}

\end{document}